\documentclass[11pt]{amsart}
\usepackage{amsfonts,amssymb,amscd,amstext,mathrsfs}
\usepackage[utf8]{inputenc}
\usepackage{hyperref}
\usepackage{verbatim}
\usepackage{graphics}
\usepackage{graphicx} 
\usepackage{times}
\usepackage{enumerate}
\usepackage[up,bf]{caption}
\usepackage{color}
\usepackage{multicol}
\input xy
\xyoption{all}

\usepackage[color=blue!20]{todonotes}

\newtheorem{theorem}{Theorem}[section]

\newtheorem{corollary}[theorem]{Corollary}

\theoremstyle{definition}
\newtheorem{definition}[theorem]{Definition}

\newtheorem{problem}[theorem]{Problem}
\newtheorem{example}[theorem]{Example}

\numberwithin{equation}{section}
\numberwithin{figure}{section}

\newcommand\Acal{\mathcal{A}}

\newcommand\Ascr{\mathscr{A}}

\newcommand\Cscr{\mathscr{C}}

\newcommand\Fscr{\mathscr{F}}
\newcommand\Gscr{\mathscr{G}}

\newcommand\Jscr{\mathscr{J}}

\newcommand\Oscr{\mathscr{O}}

\newcommand\A{\mathbb{A}}
\newcommand\B{\mathbb{B}}
\newcommand\C{\mathbb{C}}
\newcommand\D{\overline{\mathbb D}}
\newcommand\CP{\mathbb{CP}}

\renewcommand\D{\mathbb D}

\newcommand\N{\mathbb{N}}

\newcommand\R{\mathbb{R}}

\newcommand\Z{\mathbb{Z}}

\newcommand\ggot{\mathfrak{g}}

\newcommand\igot{\mathfrak{i}}

\renewcommand\igot{\mathfrak{i}}

\newcommand\Mgot{\mathfrak{M}}
\newcommand\Ngot{\mathfrak{N}}

\renewcommand\imath{\igot}

\newcommand\hra{\hookrightarrow}

\newcommand\longhookrightarrow{\ensuremath{\lhook\joinrel\relbar\joinrel\rightarrow}}

\newcommand\wh{\widehat}
\newcommand\di{\partial}
\newcommand\dibar{\overline\partial}

\newcommand\Id{\mathrm{Id}}

\newcommand\TC{\mathrm{TC}}

\def\Ell1{\mathrm{Ell_1}}
\def\DEll1{\mathrm{DEll_1}}

\numberwithin{equation}{section}

\begin{document}
\title{From Stein manifolds to Oka manifolds: \\ the h-principle in complex analysis}
\author{Franc Forstneri\v{c}}

\address{Franc Forstneri\v{c}, Faculty of Mathematics and Physics, University of Ljubljana, Jadranska 19, 1000 Ljubljana, Slovenia}

\address{Franc Forstneri\v{c}, Institute of Mathematics, Physics, and Mechanics, Jadranska 19, 1000 Ljubljana, Slovenia}

\email{franc.forstneric@fmf.uni-lj.si}

\thanks{Research supported by the 
European Union (ERC Advanced grant HPDR, 101053085) and
by grants P1-0291 and N1-0237 from ARIS, Republic of Slovenia} 


\date{24 September 2025}


\begin{abstract} This introduction to
the homotopy principle in complex analysis and geometry, 
better known as the Oka theory, is aimed at wide mathematical audience. 
After a brief historical survey of the h-principle in smooth analysis 
and geometry, I present the key notions of Oka manifolds and Oka maps, 
which developed from the Oka--Grauert principle and Gromov's theory 
of elliptic complex manifolds and elliptic holomorphic submersions. 
I discuss recent and ongoing developments, open problems, 
and mention some applications.
The paper also includes a brief survey of the recently developed 
h-principles in the classical theory of minimal surfaces.  
\end{abstract}

\maketitle

%
%
%
%
\section{Introduction.}\label{sec:intro} 

The homotopy principle, or the h-principle for short, is said to hold in an  
analytic problem if a solution exists in the absence of topological 
or homotopy theoretic obstructions.
Most instances concern solutions of underdetermined 
partial differential equations, PDEs, and more generally of 
partial differential relations, PDRs. In the holomorphic case, 
we can talk of holomorphic partial differential relations, HPDRs.
A PDE or (H)PDR is said to satisfy the h-principle if any formal
(non-holonomic) solution can be deformed 
through formal solutions to a genuine solution.
Thus in the presence of the h-principle, a differential topological or a complex analytic problem reduces to a homotopy theoretic problem. 
We speak of the basic h-principle, pertaining to individual
solutions, and of the parametric h-principle. The latter involves
a pair of parameter spaces $Q\subset P$, which are typically compact
Hausdorff spaces, and a continuous family of formal solutions parametrised
by points $p\in P$, which are genuine solutions for $p\in Q$.
The parametric h-principle is said to hold if any such family can be 
deformed to a family of genuine solutions while preserving 
the solutions for values $p\in Q$. 
This implies that the inclusion of the space 
of genuine solutions in the space of formal solutions 
is a weak homotopy equivalence, 
that is, it induces a bijection of path components and an isomorphism 
of $k$-th homotopy groups of the two spaces for every $k\in \N=\{1,2,\ldots\}$.
Sometimes there is even a genuine homotopy 
equivalence between the two spaces.
On the other hand, the failure of the h-principle means that the objects 
being studied have nontrivial geometry which cannot be reduced to 
purely topological considerations. 

The first known h-principle is the 
Whitney--Graustein theorem \cite{Whitney1937CM} which 
says that two immersions of a circle in the plane are regularly homotopic 
(connected by a path of immersions) if and only if they have the 
same winding number. Two decades later, Smale \cite{Smale1959} 
classified regular homotopy classes of immersions 
of the $k$-sphere $S^k$ in the Euclidean space 
$\R^n$ for $1\le k <n$ by elements
of the $k$-th homotopy group $\pi_{k}(V_{k}(\R^n))$ of the Stiefel manifold 
of $k$-frames in $\R^n$. It turns out that for $n\in\{2,6\}$ any two immersions
$S^n\to \R^{n+1}$ are regularly homotopic, so one can turn these 
spheres inside-out through immersions in $\R^{n+1}$. 
Hirsch \cite{Hirsch1959} proved the h-principle for 
immersions of smooth $k$-manifolds to $\R^n$ for $k<n$,
and also for $k=n$ in the case of open manifolds. 
Phillips developed the h-principle for smooth submersions 
and foliations \cite{Phillips1967,Phillips1968,Phillips1969}. 

%
%
In another direction, Nash \cite{Nash1954,Nash1956} 
and Kuiper \cite{Kuiper1955} proved that a smooth 
compact Riemannian manifold $(M,g)$ of dimension $k$ which 
admits an embedding in $\mathbb{R}^n$ for some $n>k$ also admits a 
$\mathscr{C}^1$ isometric embedding in $\mathbb{R}^n$ 
with its Euclidean metric. In fact, every short embedding can be 
uniformly approximated by isometric embeddings. This is false 
for $\mathscr{C}^2$ isometric embeddings 
due to curvature obstructions in low codimension. 

A major conceptual development of the theory was made by 
Mikhail Gromov \cite{Gromov1973,Gromov1986}, 
and the term h-principle entered into general use following his work. 
Gromov formulated the h-principle in very general terms 
and developed important new techniques 
for solving it, some of them together with Yakov Eliashberg. 
We refer to the surveys 
\cite{EliashbergMishachev2002,Gromov1986,Spring1998}
for these developments.

%
%
%
\smallskip
\noindent {\bf The h-principle in complex analysis.} 
The first instance of the h-principle in complex analysis 
is the result of Kiyoshi Oka \cite{Oka1939} from 1939, 
saying that on a domain of holomorphy in $\C^n$ 
every topological complex line bundle admits a compatible structure
of a holomorphic line bundle, and any two holomorphic lines bundles which 
are topologically isomorphic are also holomorphically isomorphic. 
A far reaching generalisation was obtained by Hans Grauert 
\cite{Grauert1957I,Grauert1957II,Grauert1958MA} in 1958. 
Recall that a {\em Stein manifold} is a complex manifold 
which admits many global holomorphic functions; 
see Sect.\ \ref{sec:Stein}. A Stein space is a complex space 
(with singularities) which enjoys the analogous function-theoretic 
properties as Stein manifolds.  
A domain in $\C^n$ is Stein if and only if it is a domain of holomorphy. 
A Riemann surface is Stein if and only if it is not compact. 
Grauert proved that complex vector bundles
and, more generally, principal fibre bundles with complex Lie group fibres
on Stein spaces have the same classification in the topological 
and the holomorphic categories.
See also the expositions by Cartan \cite{Cartan1958} and 
Henkin and Leiterer \cite{HenkinLeiterer1998,Leiterer1990}, 
and the subsequent works of Forster and Ramspott 
\cite{ForsterRamspott1966IM2,ForsterRamspott1966IM1,ForsterRamspott1968IM,ForsterRamspott1968AM}. 
These results led to the heuristic {\em Oka principle}, 
a term coined by J-P.\ Serre stating that 
{\em analytic problems on Stein spaces which can be cohomologically
formulated have only topological obstructions.}  

A homotopy-theoretic viewpoint on Oka theory was initiated 
by Gromov \cite{Gromov1986,Gromov1989} in the late 1980s. 
Modern Oka theory focuses on 
properties of a complex manifold $Y$ which imply 
that any continuous map $X\to Y$ from a Stein space 
$X$ is homotopic to a holomorphic map; the same  
problem is considered for sections of holomorphic maps 
onto Stein spaces. (All Stein spaces will be assumed reduced.) 
Adding the approximation and interpolation conditions  
modelled on properties of holomorphic functions on Stein spaces, 
one obtains several Oka properties (see Def.\ \ref{def:OP})
which a given complex manifold may or may not have. 
One of Gromov's main results in \cite{Gromov1989} 
is that a geometric condition called {\em ellipticity} --- the existence of a 
dominating holomorphic spray on the manifold, see Def.\ \ref{def:elliptic} --- 
implies all these Oka properties.
The analogous result holds for sections of elliptic submersions 
onto Stein spaces; see \cite{Gromov1989,ForstnericPrezelj2002}  
and \cite[Theorem 6.2.2]{Forstneric2017E}.
The Oka--Grauert theory fits in this framework by considering maps 
to classifying spaces, which are complex homogeneous 
and hence elliptic. However, Gromov's theorem, with 
detailed presentations in \cite{ForstnericPrezelj2000,ForstnericPrezelj2002,Forstneric2010PAMQ} 
and \cite[Chaps.\ 5--6]{Forstneric2017E}, 
applies to a considerably wider class of problems.

Subsequent research was motivated by the question whether
ellipticity is also a necessary condition for the Oka principle,
in finding new examples of elliptic manifolds and 
submersions, and in attempts to unify the theory by finding  
simpler conditions characterising the Oka property.
The first question remained open until recently and we shall
say more about it Sect.\ \ref{sec:elliptic}, where we also mention 
examples of elliptic manifolds. 
The last question was resolved by 
the introduction of the following Runge-type 
approximation property in \cite{Forstneric2006AM}.

%
%
\begin{definition}\label{def:CAP}
A complex manifold $Y$ enjoys the {\em Convex Approximation Property, CAP,} 
if every holomorphic map from an open neighborhood of 
a compact convex set $K\subset \C^n$, $n\in \N$, 
to $Y$ can be approximated uniformly on $K$ by holomorphic 
maps $\C^n\to Y$.
\end{definition}

The following is a cumulative result of 
\cite{Forstneric2005AIF,Forstneric2006AM,Forstneric2009CR,Forstneric2010PAMQ}; see also \cite[Theorem 5.4.4 and Prop.\  5.15.1]{Forstneric2017E}. 

%
%
\begin{theorem}\label{th:Oka}
A complex manifold $Y$ which enjoys CAP has all Oka properties in 
Def.\ \ref{def:OP} for maps $X\to Y$ from any Stein space $X$.
Furthermore, all these Oka properties are pairwise equivalent.
\end{theorem}

A complex manifold $Y$ satisfying these 
equivalent conditions is called an {\em Oka manifold} 
\cite{Forstneric2009CR,Larusson2004}. 
This gives an affirmative answer to Gromov's problem 
\cite[3.4 (D)]{Gromov1989}, asking whether the Oka property can be  
characterised by a simple Runge-type approximation property. 
Every elliptic manifold in the sense of 
Gromov (see Def.\ \ref{def:elliptic}) enjoys CAP  
and hence is an Oka manifold; see Theorem \ref{th:ellCAP}. 
This extends Grauert's result that every complex homogeneous 
manifold is Oka.  Essentially the same proof shows the following;
see \cite[Theorem 5.4.4]{Forstneric2017E}. 

\begin{theorem}\label{th:Okafibre}
If $h:Z\to X$ is a stratified holomorphic fibre bundle over a Stein space $X$
whose fibres are Oka manifolds, then sections $X\to Z$ of $h$ 
satisfy all forms of the Oka principle.
\end{theorem}

It is also important to consider the Oka principle
for lifting holomorphic maps in the following diagram.
\[
    \xymatrix{  & Y \ar[d]^{h} \\ X \ar[r]^{\ \ f} \ar@{-->}[ur]^{F} & Z}
\]
A holomorphic map $h:Y\to Z$ of complex manifolds
is said to be an {\em Oka map} if it is a topological (Serre) fibration
and it satisfies the following Oka property. 

\smallskip
\noindent {\em The Oka property for liftings}:
 Given a holomorphic map $f:X\to Z$ from a Stein space $X$, 
any continuous lifting $X\to Y$ of $f$ is homotopic  
to a holomorphic lifting $F:X\to Y$, with approximation and
interpolation conditions in Def.\ \ref{def:OP}, and the corresponding 
result holds in the parametric setting. 
\smallskip

We refer to \cite[Definition 7.4.1]{Forstneric2017E} for the precise definition.
A holomorphic fibre bundle with Oka fibre is an Oka map
(see \cite[Corollary 7.4.8]{Forstneric2017E}). 
In particular, a holomorphic covering map is an Oka map. 
The importance of Oka maps lies in 
the following result; see \cite[Theorem 3.15]{Forstneric2023Indag}.

%
%
\begin{theorem}\label{th:Okamap}
If $Y\to Z$ is an Oka map of connected complex manifolds, then 
$Y$ is an Oka manifold if and only if $Z$ is an Oka manifold.
\end{theorem}

L\'arusson \cite{Larusson2004,Larusson2005} 
constructed a model category in which Stein manifolds are cofibrant,
Stein inclusions are cofibrations, Oka manifolds are fibrant, 
and Oka maps are fibrations. 
It is thus fair to say that Stein manifolds are the natural sources 
of holomorphic maps while Oka manifolds are their natural targets.

Recently, Kusakabe showed in \cite[Theorem 1.3]{Kusakabe2021IUMJ}
that the following property also characterises Oka manifolds. 

%
%
\begin{definition} \label{def:CRE}
A complex manifold $Y$ enjoys {\em Convex Relative Ellipticity, CRE}, 
if for every compact convex set $K\subset\C^n$, $n\in \N$, and
holomorphic map $f:U' \to Y$ from an open neighbourhood of $K$
there exist a neighbourhood $U\subset U'$ of $K$
and a holomorphic map $F:U\times \C^N\to Y$ for some $N\in \N$ 
such that for every $z\in U$ we have $F(z,0)=f(z)$ and the differential
$dF(z,\cdotp)|_0: T_0\C^N\to T_{f(z)}Y$ at $0\in\C^N$ is surjective.
\end{definition}

A map $F$ as in the above definition is called {\em a dominating holomorphic
spray with the core} $F(\cdotp,0)=f$. This condition, for maps from 
any Stein manifold in place of a convex domain $U\subset \C^n$,
was introduced by Gromov \cite{Gromov1989} as condition $\Ell1$
(see also \cite[Definition 3.1]{Forstneric2023Indag}).
It has also been called {\em relative ellipticity} in the literature to emphasise
that it concerns sprays of holomorphic maps to the given manifold,
as opposed to sprays on the manifold itself.
It is easily seen that every holomorphic map $f:U'\to Y$ in Def.\ \ref{def:CRE}
admits a local dominating 
spray $F:U\times \B\to Y$ with the core $f|_U$, where 
$K\subset U\Subset U'$, $U$ is convex, and $\B$ is a ball in some $\C^N$. 
The main point is to find a global spray $U\times\C^N\to Y$ which 
agrees with $F$ to the second order along $U\times \{0\}$.
This shows that relative ellipticity, and hence CRE,
is implied by the basic Oka property with jet interpolation (BOPJI); 
see Def.\ \ref{def:OP} (iv'). By definition, every Oka manifold satisfies BOPJI 
and hence also CRE. To complete the picture, 
Kusakabe proved in \cite{Kusakabe2021IUMJ} that CRE implies CAP.  
(An exposition of his proof can also be found in 
\cite[Theorem 3.3]{Forstneric2023Indag}.) 
Summarising, the following properties of a complex manifold are equivalent.

\begin{theorem}\label{th:OkamainI}
$
	\mathrm{OKA} \ \Longleftrightarrow\  \mathrm{CAP} 
	\Longleftrightarrow\  \mathrm{CRE}. 
$
\end{theorem}

The characterisation of Oka manifolds
by CRE is sometimes easier to apply then CAP.
By using it, Kusakabe showed in 
\cite[Theorem 1.4]{Kusakabe2021IUMJ}
that the Oka property is Zariski local in the following sense.

%
%
\begin{theorem}
\label{th:localisation}
If every point in a complex manifold $Y$ admits a 
Zariski open$\,$\footnote{A domain in a complex manifold $Y$ is 
Zariski open if its complement is a closed complex subvariety of $Y$.} 
Oka neighbourhood, then $Y$ is an Oka manifold.
\end{theorem}

This result has been an important source of 
new examples of Oka manifolds in recent times.
In \cite{Kusakabe2021MZ}, Kusakabe also characterised the Oka property 
of sections of holomorphic submersions by a fibred version of CRE.

Oka manifolds have a number of special analytic properties. 
It is immediate that an Oka manifold $Y$ is dominable by the Euclidean 
space $\C^n$ with $n={\dim Y}$ at every point $y\in Y$, 
in the sense that there exists
a holomorphic map $f:\C^n \to Y$ with $f(0)=y$ 
such that $df_0:T_0\C^n\to T_yY$ is an isomorphism.
Conversely, dominability at most points implies the Oka property for maps 
from open Riemann surfaces (see Theorem \ref{th:Oka1}), 
but it is not known whether it implies the full Oka property. 
With some more work, one can see that a connected Oka manifold $Y$ 
admits a surjective holomorphic map $f:\C^n \to Y$, $n=\dim Y$,  
such that for every $y\in Y$ there is a point $z\in f^{-1}(y)\in \C^n$ 
such that $df_z:T_z\C^n\to T_yY$ is an isomorphism 
(see \cite[Theorem 1.1]{Forstneric2017Indam}). It follows that 
a connected Oka manifold does not admit any nonconstant 
bounded plurisubharmonic function. Furthermore, 
all pseudometrics on complex manifolds which are
decreasing under holomorphic maps and vanish on Euclidean spaces
(such as the Kobayashi metric, the Caratheodory metric, 
the Eisenmann metrics, the Sibony metric, to name some of 
the best known ones) also vanish on any Oka manifold. 
Since a compact complex manifold of general type 
(i.e.\ of maximal Kodaira dimension equal to its complex 
dimension) is not dominable by a Euclidean space 
(see Kobayashi and Ochiai \cite{KobayashiOchiai1975}), 
no such manifold is Oka. 
There are many examples of compact Oka manifolds among 
Fano manifolds (Kodaira dimension $-\infty$), 
although it is not known whether every Fano manifold is Oka.
We refer to \cite[Sect.\ 7.3]{Forstneric2017E} for a review of 
known Oka manifolds among compact complex surfaces. 
For Riemann surfaces there is a precise dichotomy --- either the surface is Oka, and then it is one of the surfaces $\CP^1$, $\C$, $\C^*=\C\setminus \{0\}$ or a torus, or it is Kobayashi hyperbolic and hence
a quotient of the disc $\D=\{z\in\C:|z|<1\}$.  

Elliptic manifolds (see Def.\ \ref{def:elliptic})
remain an important source of examples
of Oka manifolds, especially among the algebraic ones. 
The theory of algebraically elliptic manifolds 
is an active area of contemporary research; see 
Sect.\ \ref{sec:elliptic}. We have recently 
shown with L\'arusson \cite{ForstnericLarusson2025MRL} 
that every projective Oka manifold is elliptic; see Theorem \ref{th:projective}.
This gives an affirmative answer to 
a conjecture of Gromov \cite[3.2.A']{Gromov1989}.
On the other hand, due to Kusakabe
\cite{Kusakabe2024AM} we now have many examples of 
noncompact Oka manifolds which are not elliptic;
see Corollary \ref{cor:complements}.
It remains an open question whether there exist compact 
nonprojective Oka manifolds which fail to be elliptic.

There are several open problems concerning 
the relationship between the Oka property and topological, analytical, and 
metrical properties of complex manifolds.
A Stein manifold $X$ has the homotopy type
of a CW complex of dimension at most $\dim_\C X=\frac12 \dim_\R X$.
Conversely, given a smooth almost complex manifold $(X,J_0)$ of 
real dimension $2n\ne 4$ with the homotopy type of a CW 
complex of dimension $\le n$, 
$J_0$ is homotopic to a Stein structure $J$ on $X$; 
see Eliashberg \cite{Eliashberg1990}.
The analogous result holds if $\dim_\R X = 4$ after a suitable change
of the $\Cscr^\infty$ structure on $X$; 
see Gompf \cite{Gompf1998,Gompf2005}. 
Furthermore, every continuous map
$f:X\to Y$ from a Stein manifold $X$ to an arbitrary complex
manifold $Y$ is homotopic to a holomorphic map with respect
to a Stein structure $J$ on $X$ which is homotopic to the original one
(and up to a change of the $\Cscr^\infty$ structure on $X$ if
$\dim_\R X=4$); see \cite{ForstnericSlapar2007MZ,ForstnericSlapar2007MRL}.
This is called the {\em soft Oka principle}. See
also the monograph by J\"oricke \cite{Joricke2025} which studies
this phenomenon in the context of braids. 
On the other hand, almost nothing is known concerning the 
homotopy types of Oka manifolds. In particular, we do not know 
which groups arise as fundamental groups of Oka manifolds.
Another intriguing question concerns the geometric
shape of Oka domains in complex manifolds; 
see Problem \ref{prob:geometric}. 

A challenging problem is to understand the relationship between
the Oka property and metrical properties of 
complete hermitian and K\"ahler manifolds. 
Since negativity of the holomorphic sectional
curvature is related to Kobayashi hyperbolicity (see  
\cite{GrauertReckziegel1965,GreeneWu1979,BroderIliashenkoMadnick2025}), 
which is the very opposite of the Oka property, 
one might expect that positivity together with completeness
is related to the Oka property. 
An important result in this direction is given by Mok's solution 
of the generalised Frankel conjecture \cite{Mok1988},
which implies that every compact {K}\"{a}hler manifold 
of nonnegative holomorphic bisectional curvature is an Oka
manifold \cite[Theorem 11.4]{Forstneric2023Indag}.  
Our recent result with Kusakabe \cite{ForstnericKusakabe2025} 
on the Oka property of disc tubes in semipositive 
ample line bundles on flag manifolds 
(see Theorem \ref{th:density2}) also contributes to the heuristic principle 
that metric positivity is related to the Oka property. 

Another problem is the relationship
between Oka manifolds and special manifolds in the sense of Campana 
\cite{Campana2004AIF,Campana2004AIF-2}. It was shown by 
Campana and Winkelmann \cite{CampanaWinkelmann2015}
that every projective manifold $Y$ which satisfies the basic 
Oka property (every continuous map from a Stein manifold to $Y$ 
is homotopic to a holomorphic map) is special.
In particular, every projective Oka manifold is special
but the converse is an open problem.

There is a growing list of applications of Oka theory, some
of which are discussed in my book \cite{Forstneric2017E}.
There are the dimensionwise optimal embedding and immersion
theorems for Stein manifolds into Euclidean spaces, 
due to Eliashberg and Gromov
\cite{EliashbergGromov1992} and Sch\"urmann 
\cite{Schurmann1997}, 
the h-principle for holomorphic immersions of Stein manifolds
$X^n\to \C^N$ for $N>n$ (see Eliashberg and Gromov 
\cite{EliashbergGromov1971} and \cite[Sect.\ 2.1.5]{Gromov1986}), 
and the h-principle for holomorphic submersions $X^n\to \C^q$
with $n>q\ge 1$ (see Forstneri\v c \cite{Forstneric2003AM}).
There exist proper holomorphic maps, immersions  
and embeddings $X\to \C^n$ of Stein manifolds for suitable 
values of $n>\dim X$ whose images avoid any given unbounded 
strictly convex set \cite{DrinovecForstneric2023JMAAproper}
and have small limit sets at infinity \cite{Forstneric2024PAMS}. 
Kutzschebauch, L{\'a}russon, and Schwarz developed  
$G$-equivariant Oka theory and the notion of a 
$G$-Oka manifold, where $G$ is a complex Lie group; see 
their survey \cite{KutzschebauchLarussonSchwarz2022CAS}.
There are major recent developments on 
factorisation problems for 
holomorphic matrix-valued functions on Stein spaces 
\cite{IvarssonKutzschebauch2012AM,IvarssonKutzschebauchLow2020,HuangKutzschebauchSchott2024,HuangKutzschebauchTran2025},
which rely on the Oka principle for sections 
of stratified elliptic submersions onto Stein spaces; 
see \cite[Theorem 6.2.2]{Forstneric2017E}.
Oka theory combined with methods of convex integration theory was used to 
establish h-principles in the classical 
theory of minimal surfaces in Euclidean spaces (see 
\cite{AlarconForstneric2014IM,ForstnericLarusson2019CAG,AlarconForstnericLarusson2024,AlarconLarusson2024complete,AlarconLarusson2025Crelle} 
and the monograph \cite{AlarconForstnericLopez2021}),   
and in the theory of holomorphic Legendrian curves 
\cite{ForstnericLarusson2018MZ,ForstnericLarussonIUMJ2022}.
A brief discussion of the h-principles in the theory of minimal surfaces
is included in Sect.\ \ref{sec:minimal}.

Another topic intimately related to Oka theory 
is the Anders\'en--Lempert theory and Varolin's density 
property concerning Stein manifolds with a large group 
of holomorphic automorphisms; see Subsec.\ \ref{ss:complements}. 
Their main feature is the approximation theorem for isotopies of 
biholomorphic maps between pseudoconvex
Runge domains by isotopies of holomorphic automorphisms
(see Forstneri\v c and Rosay \cite[Theorem 1.1]{ForstnericRosay1993}, 
\cite[Theorem 4.10.5]{Forstneric2017E}). 
Surprisingly, such highly symmetric objects exist in 
big classes and were found by applying 
criteria developed by Kaliman, Kutzschebauch and coauthors 
over the last two decades; see the surveys 
\cite{ForstnericKutzschebauch2022,Kutzschebauch2020} 
and \cite[Chap.\ 4]{Forstneric2017E}. This theory had a great impact 
on the solution of famous old problems, among them the 
holomorphic linearisation problem, the embedding problem for Riemann surfaces, 
the h-principle for holomorphic submersions of Stein manifolds 
to Euclidean spaces, and many others. A Stein manifold $Y$ 
with the density property is Oka at infinity 
(see Theorem \ref{th:complements} due to Kusakabe \cite{Kusakabe2024AM}), 
and it contains properly embedded copies 
of any Stein manifold $X$ with $2\dim X<\dim Y$ in every homotopy class 
of maps $X\to Y$ \cite{AndristForstnericRitterWold2016}.
Arosio and L\'arusson studied the dynamics of 
generic endomorphisms of Oka--Stein manifolds \cite{ArosioLarusson2022}
and of Stein manifolds with the density property 
\cite{ArosioLarusson2023,ArosioLarusson2025}. 

Reader, follow me to the continuation of this story. 

%
%
\section{Stein manifolds.}\label{sec:Stein}
In this section we recall the basic properties of Stein manifolds
which are relevant for this survey, referring to
\cite{GunningRossi2009,Hormander1990,GrauertRemmert1979,Forstneric2017E} 
for more information. Stein manifolds were introduced in the literature 
by Karl Stein in 1951 \cite{Stein1951} under the name 
{\em holomorphically complete manifolds} by the following axioms. 

\vspace{-1mm}
\begin{enumerate}[\rm (a)]
\item Holomorphic functions on $X$ separate any distinct pair of points.
\vspace{-1mm}
\item For every point $p\in X$ there is a holomorphic map
$f=(f_1,\ldots,f_n):X\to\C^n$ with $n=\dim_\C X$ whose 
differential at $p$ has complex rank $n$.
\vspace{-1mm}
\item For every compact set $K$ in $X$, its 
holomorphically convex (or $\Oscr(X)$-convex) hull 
\begin{equation}\label{eq:hull}
	\wh K=\{p\in X: |f(p)|\le \max_{x\in K}|f(x)|\ \text{for all}\ f \in \Oscr(X)\}
\end{equation}
is also compact. (By $\Oscr(X)$ we denote the space of holomorphic functions $X\to\C$.)
\end{enumerate}

A compact set $K\subset X$ is said to be holomorphically convex in $X$,
or $\Oscr(X)$-convex, if and only if $K=\wh K$. 
When $X=\C^n$, such a set is called polynomially convex.
The same axioms define Stein spaces if (b) is replaced by the condition 
that holomorphic functions on $X$ generate the ring 
$\Oscr_{X,x}$ of germs of holomorphic functions at every point $x\in X$. 
In fact, (b) is a consequence of (a) and (c).
A complex space satisfying condition (c) is said to be holomorphically convex. 
A domain in $\C^n$ is Stein if and only if it is holomorphically convex 
if and only if it is a domain of holomorphy. 
An open Riemann surface is Stein according to Behnke and Stein
\cite{BehnkeStein1949}.

Stein manifolds have a number of characterisations. 
A connected complex manifold $X$ is Stein if and only if it
is biholomorphic to a closed complex submanifold of a 
Euclidean space $\C^N$;
we can take $N=2\dim X+1$. Equivalently,
$X$ admits a proper holomorphic embedding $X\hra\C^{2n+1}$.
See \cite[Theorem 2.4.1]{Forstneric2017E} for a summary 
and references. Thus, Stein manifolds are complex analytic analogues 
of affine algebraic manifolds. 
A more precise embedding theorem of Eliashberg and Gromov
\cite{EliashbergGromov1992} and Sch\"urmann 
\cite{Schurmann1997} is that a Stein manifold of dimension $n>1$ 
embeds properly holomorphically in $\C^N$ with 
$N=\left[\frac{3n}{2}\right]+1$; 
this dimension is optimal by examples of Forster \cite{Forster1970}.  
Similar embedding results hold for Stein spaces with bounded
local embedding dimension; see \cite[Theorem 9.3.7]{Forstneric2017E}.
The proofs rely on Oka theory; see the detailed exposition 
in \cite[Secs.\ 9.3--9.4]{Forstneric2017E}. 
The problem whether every open Riemann surface embeds
as a closed nonsingular complex curve in $\C^2$ is still wide open;
see the summary of known results
in \cite[Secs.\ 9.10--9.11]{Forstneric2017E} and the recent paper
\cite{DiSalvoRitterWold2023}. 

Another important characterisation is that a complex manifold $X$
is Stein if and only if it admits a smooth strongly plurisubharmonic
exhaustion function $\rho:X\to[0,+\infty)$, 
i.e., one satisfying $dd^c\rho>0$; see Grauert \cite{Grauert1958AM}.
The analogous condition characterises Stein spaces 
with singularities, see Forn\ae ss and Narasimhan 
\cite{FornaessNarasimhan1980}. 

One of the highlights of complex analysis on Stein spaces
are Cartan's Theorems $A$ and $B$.
Let $\Fscr$ be a coherent analytic sheaf on a Stein space $X$.
Theorem A says that every stalk $\Fscr_x$ for $x\in X$ is generated
as an $\Oscr_{X,x}$-module by germs of global 
sections of $\Fscr$, with Runge approximation of sections 
on compact $\Oscr_X$-convex subsets of $X$. 
Theorem B says that $H^q(X,\Fscr)=0$
for all $q=1,2,\ldots$. Furthermore, 
the Dolbeault cohomology groups of a Stein manifold 
vanish, that is, the Cauchy--Riemann equation $\dibar u=\alpha$ is solvable 
on $X$ for any $(p,q)$-form $\alpha$ $(p\ge 0,\ q\ge 1)$
with $\dibar \alpha=0$. 

An important role in Oka theory is played by the 
theorem of Siu \cite{Siu1976} saying that 
a locally closed Stein subvariety $X$ in a complex space $Z$
admits an open Stein neighbourhood. 
See also \cite[Secs.\ 3.1--3.2]{Forstneric2017E} and 
the generalisation of Siu's theorem in \cite[Theorem 3.2.1]{Forstneric2017E}. 
A related result of Poletsky \cite{Poletsky2013}
concerns Stein neighbourhoods of graphs of holomorphic mappings
that are continuous up to the boundary of a domain.

The following classical results on holomorphic functions on a 
Stein space $X$ motivate the definition of Oka properties
and of Oka manifolds; see Def.\ \ref{def:OP} and Theorem \ref{th:Okamain}.
\begin{enumerate}[\rm (i)]
\item {\em The Oka--Weil approximation theorem}: 
If $K$ is a compact $\Oscr(X)$-convex subset of $X$ then any
holomorphic function $U\to \C$ on a neighbourhood $U\subset X$ of $K$
is a uniform limit on $K$ of holomorphic functions $X\to\C$.
\vspace{-3mm}
\item {\em The Cartan--Oka extension theorem:} 
If $X'$ is a closed complex subvariety of $X$ then any holomorphic
function $X'\to \C$ extends to a holomorphic function $X\to\C$.
\vspace{-1mm}
\item If $K\subset U\subset X$ and $X'\subset X$ are as above and
$f:U\cup X'\to \C$ is a continuous function which is 
holomorphic on $U$ and on $X'$, we can 
approximate $f$ uniformly on $K$ by 
functions $F\in \Oscr(X)$ satisfying $F|_{X'}=f|_{X'}$. 
\vspace{-1mm}
\item If a function $f$ in (iii) is holomorphic on a neighbourhood 
of $K\cup X'$, then $F\in \Oscr(X)$ in (iii) can be chosen to approximate
$f$ on $K$ and to agree with $f$ to any given finite order along $X'$.  
\end{enumerate}
Note that (iii) is a combination of (i) and (ii), while (iv) upgrades the
interpolation on $X'$ to jet interpolation. The analogous results
hold in the parametric case; see \cite[Theorem 2.8.4]{Forstneric2017E}.

%
%
%
%
\section{Oka manifolds.}\label{sec:Oka}
Results (i)--(iv) above, and their parametric 
analogues in \cite[Theorem 2.8.4]{Forstneric2017E}, 
describe the essential function-theoretic properties of Stein spaces.
We now consider the analogous conditions 
for maps $X\to Y$ from a Stein space $X$ to a complex manifold $Y$.

\subsection{Oka properties characterising Oka manifolds.}

%
%
\begin{definition}\label{def:OP} 
Let $Y$ be a complex manifold.
\vspace{-1mm}
\begin{enumerate}[\rm (i')]
\item $Y$ has the {\em basic Oka property with approximation, BOPA}, 
if for any Stein space $X$, compact holomorphically
convex subset $K\subset X$, and continuous map 
$f:X\to Y$ which is holomorphic on a neighbourhood of $K$ 
there is a homotopy $f_t:X\to Y$ $(t\in I=[0,1])$ such that 
$f_0=f$, each map $f_t$ is holomorphic on a neighbourhood of 
$K$ (independent of $t\in I$) and approximates $f$ uniformly on $K$ 
and uniformly in $t\in I$ to a given precision, 
and the map $f_1$ is holomorphic on $X$.
\vspace{-1mm}
\item 
$Y$ has the {\em basic Oka property with interpolation, BOPI}, 
if for any Stein space $X$, closed complex subvariety $X'\subset X$, 
and continuous map $f:X\to Y$ such that $f|_{X'}$ is holomorphic 
there is a homotopy $f_t:X\to Y$ $(t\in I)$ such that 
$f_0=f$, $f_t|_{X'}=f|_{X'}$ for every $t\in I$, and $f_1$ is holomorphic on $X$.
\vspace{-1mm}
\item 
$Y$ has the
{\em basic Oka property with approximation and interpolation, BOPAI}, 
if the combination of BOPA and BOPI holds 
(cf.\ property (iii) in Sect.\ \ref{sec:Stein}).
\vspace{-1mm}
\item 
$Y$ has the {\em basic Oka property with approximation 
and jet interpolation, BOPAJI}, if the analogue of property (iv) in 
Sect.\ \ref{sec:Stein} holds for maps $X\to Y$. Omitting 
the approximation condition defines the 
{\em basic Oka property with jet interpolation, BOPJI}.
\end{enumerate}
\end{definition}

In order to avoid topological obstructions, the initial continuous map $f:X\to Y$ (or a continuous family of maps in the parametric case) is assumed to be defined on all of $X$. Note that CAP (see Def.\ \ref{def:CAP}) is a restricted version of BOPA, 
and CRE (see Def.\ \ref{def:CRE}) is a restricted version of BOPJI. 
Similarly one defines the parametric versions of these properties 
called POPA, POPI, POPAI, POPAJI and POPJI, respectively;
see \cite[Sect.\ 5.15]{Forstneric2017E}.
This gives several ostensibly different Oka-type conditions on a complex manifold. It turns out that they are pairwise equivalent 
\cite[Proposition 5.15.1]{Forstneric2017E}, although several implications
are highly nontrivial. 

%
%
\begin{theorem}
\label{th:Okamain}
Every Oka property mentioned above 
is equivalent to the Convex Approximation Property 
(CAP) in Def.\ \ref{def:CAP}, and to Convex Relative Ellipticity (CRE) 
in Def.\ \ref{def:CRE}. 
\end{theorem}

A complex manifold $Y$ is said to be an {\em Oka manifold} if and only if 
it satisfies any and hence all of these conditions.
Oka manifolds appear in Gromov's paper under the name
Ell$_\infty$ manifolds \cite[3.3.C']{Gromov1989}. 
However, the mentioned equivalences were proved much later,
and the notion of an Oka manifold was introduced in \cite{Forstneric2009CR}. 

The following is a well known consequence of the parametric h-principle. 

\begin{corollary}\label{cor:WHE}
If $X$ is a Stein space and $Y$ is an Oka manifold then the inclusion 
$
	\Oscr(X,Y) \longhookrightarrow \Cscr(X,Y)
$
of the space of holomorphic maps $X\to Y$ in the space 
of continuous maps (both endowed with the compact-open topology)
is a weak homotopy equivalence.
\end{corollary}

L\'arusson showed in \cite{Larusson2015PAMS} 
that the above inclusion is a genuine homotopy
equivalence if $X$ is a Stein manifold which admits a strongly 
plurisubharmonic exhaustion function $\rho:X\to [0,\infty)$
with only finitely many critical points.
In this case, the space $\Oscr(X,Y)$ 
is a deformation retract of $\Cscr(X,Y)$ if $Y$ is Oka.
This holds in particular if $X$ is a strongly pseudoconvex domain 
in another Stein manifold, or if it is an affine algebraic manifold.

More precise statements 
can be found in \cite[Theorem 5.4.4 and Proposition 5.15.1]{Forstneric2017E}.
The proof of Theorem \ref{th:Okamain} (without considering CRE) 
takes up \cite[Secs.\ 5.7--5.13]{Forstneric2017E}. 
The implication CRE $\Longrightarrow$ CAP is due to 
Kusakabe \cite[Theorem 1.3]{Kusakabe2021IUMJ}; 
see also the exposition in \cite[Theorem 3.3]{Forstneric2023Indag}. 

The following result (see \cite[Theorem 1.3]{Forstneric2023Indag}) gives
a more precise control of the image of a holomorphic map.
It is particularly useful in conjunction with Theorem \ref{th:complements}
and is often applied in the constructions of proper holomorphic maps. 
It is obtained by following \cite[proof of Theorem 5.4.4]{Forstneric2017E}. 
We only state the basic case.

%
%
\begin{theorem}\label{th:Oka2}
Assume that $X$ is a Stein space, 
$K$ is a compact $\Oscr(X)$-convex set in $X$, 
$X'$ is closed complex subvariety of $X$, 
$\Omega$ is an Oka domain in a complex manifold $Y$, 
and $f:X\to Y$ is a continuous map which is holomorphic 
on a neighbourhood of $K$ and on $X'$ such that 
$
	f(\overline{X\setminus K}) \subset \Omega.
$  
Then there is a homotopy $f_t:X\to Y$ $(t\in [0,1])$ connecting $f=f_0$ 
to a holomorphic map $f_1:X\to Y$ 
satisfying the conclusion of Theorem \ref{th:Okamain} 
such that $f_t(\overline{X\setminus K}) \subset \Omega$
holds for all $t\in[0,1]$.
\end{theorem}

%
%
Note that every Oka property in Def.\ \ref{def:OP} includes
an approximation or an interpolation condition. 
A complex manifold $Y$ is said to satisfy the 
{\em basic Oka property}, BOP, if every continuous map $X\to Y$ 
from a Stein manifold $X$ is homotopic to a holomorphic map. 
Every Oka manifold clearly satisfies BOP but the converse fails in general.
Indeed, any holomorphically contractible complex manifold, 
such as the disc, satisfies BOP but need not be Oka. 
Campana and Winkelmann proved in \cite{CampanaWinkelmann2015}
that a projective manifold $Y$ satisfying BOP 
is Campana special \cite{Campana2004AIF,Campana2004AIF-2}.

\begin{problem} Is every projective manifold satisfying BOP an Oka manifold? 
\end{problem}

A survey of examples of Oka manifolds which predates their characterisation 
by CRE can be found in \cite[Chaps.\ 5--7]{Forstneric2017E}.
The recent survey \cite{Forstneric2023Indag} includes
many new examples based on the equivalence 
CRE $\Longleftrightarrow$ OKA, 
Theorem \ref{th:localisation}, and the use of the density property. 
We present some of them in the following subsection.

\subsection{Oka properties of complements of holomorphically convex sets.} 
\label{ss:complements}
%
%
A complex manifold $Y$ is said to enjoy the {\em density property}
(see Varolin \cite{Varolin2000,Varolin2001}) if every holomorphic 
vector field on $Y$ can be approximated uniformly 
on compacts by Lie combinations 
of $\C$-complete holomorphic vector fields.
The Euclidean spaces $\C^n$, $n>1$, enjoy this property by 
Anders\'en and Lempert \cite{AndersenLempert1992}. One of the main
results in this subject is that if $Y$ is a Stein manifold 
with the density property then every isotopy of biholomorphic maps
$\Phi_t:\Omega\to \Phi_t(\Omega)$, $t\in [0,1]$, with $\Phi_0=\Id_\Omega$  
between pseudoconvex Runge domains in $Y$ 
can be approximated by isotopies of holomorphic automorphisms of $Y$ 
(see Forstneri\v c and Rosay
\cite{ForstnericRosay1993}, \cite[Theorem 4.10.5]{Forstneric2017E}). 
Every Stein manifold with the density property is an Oka manifold 
(see \cite[Theorem 4]{KalimanKutzschebauch2008MZ}).
Surveys can be found in \cite[Chap.\ 4]{Forstneric2017E},   
\cite{ForstnericKutzschebauch2022}, and \cite{Kutzschebauch2020}. 

In his seminal paper \cite{Kusakabe2024AM}, Kusakabe applied the 
characterisation of Oka manifolds by CRE 
(see Theorem \ref{th:Okamain}) to establish the Oka property of 
complements of certain closed holomorphically convex 
sets in Stein manifolds with the density property. This gives  
a plethora of new examples of Oka manifolds, most of which
are not Gromov elliptic. The following result combines
Theorem 1.2 and Corollary 1.3 in \cite{Kusakabe2024AM}.
(See also \cite{ForstnericWold2020MRL}.) 

%
%
\begin{theorem}\label{th:complements}
If $K$ is a compact polynomially convex set in $\C^n$, $n>1$,  
then $\C^n\setminus K$ is an Oka manifold. 
More generally, if $Y$ is a Stein manifold with the density property 
and $K\subset Y$ is a compact $\Oscr(Y)$-convex set then 
$Y\setminus K$ is an Oka manifold. 
\end{theorem}

The idea of proof is that, given a set $K\subset Y$ as in the theorem 
and a holomorphic map $f:U'\to Y\setminus K$ 
from a neighbourhood of a compact convex
set $L\subset\C^N$, one can find a neighbourhood
$U\subset U'$ of $L$ and a holomorphic map
$F:U\times \C^n \to Y\setminus K$, with $n=\dim Y$, such that
for every $z\in U$ the map
$F_z=F(z,\cdotp):\C^n \to Y$ is biholomorphic onto its
image in $Y\setminus K$ and satisfies $F_z(0)=f(z)$ 
(see \cite[Theorems 1.1 and 3.1]{ForstnericWold2020MRL}).  
Thus, $\{F_z(\C^n)\}_{z\in U}$ is a holomorphically varying family of 
Fatou--Bieberbach domains in $Y\setminus K$.
Clearly, such $F$ is a dominating spray with the core $f$, 
so $Y\setminus K$ satisfies CRE and hence is an Oka manifold. 
The construction of $F$ uses the density property of $Y$ and the 
aforementioned approximation of isotopies of biholomorphic maps between 
Stein Runge domains in $Y$ by holomorphic automorphisms of $Y$
(see \cite[Theorem 4.10.5]{Forstneric2017E}). 

A closed unbounded subset $S$ of $\C^n$ is said to be polynomially convex 
if it is exhausted by an increasing sequence of compact polynomially 
convex sets. For $Y=\C^n$, $n>2$, Theorem \ref{th:complements} 
is a special case of the following result of Kusakabe 
\cite[Theorems 1.6 and 4.2]{Kusakabe2024AM}.

%
%
\begin{theorem}\label{th:unboundedsets}
If $n>2,\ c>0$ and 
$S\subset \left\{(z,w)\in \C^{n-2}\times \C^{2}: |w| \le c(1+|z|)\right\}$
is a closed polynomially convex set, then $\C^n\setminus S$ is an Oka manifold.
\end{theorem}

In \cite{ForstnericWold2024IMRN}, these techniques were used to prove that 
most concave domains in $\C^n$ for $n>1$ are Oka. 
In particular, the following holds; 
see \cite[Theorem 1.8]{ForstnericWold2024IMRN}. 

%
%
\begin{theorem} \label{th:convexnoline}
If $E$ is a closed convex set in $\C^n$, $n>1$, which does not 
contain any affine real line, then $\mathbb C^n\setminus E$ is an Oka domain. 
\end{theorem}

In particular, if $\phi:\C^{n-1}\times\R\to \R_+$ is a strictly convex
function then the concave domain
\begin{equation}\label{eq:model}
	\{z=(z',z_n) \in \C^n: \Im z_n<\phi(z',\Re z_n)\}
\end{equation}
below its graph is an Oka domain while the domain above the graph is
Kobayashi hyperbolic. Note that these domains can be 
arbitrarily close to a halfspace, the latter being neither Oka nor hyperbolic. 

Consider $\C^n$ as an affine domain in the projective space 
$\CP^n$. Given a subset $E\subset\C^n$, 
we denote by $\overline E$ its topological closure in $\CP^n$. 
Theorem \ref{th:convexnoline} is a special case of the following 
result \cite[Theorem 1.1]{ForstnericWold2024IMRN}.

%
%
\begin{theorem}\label{th:FWmain}
If $E$ is a closed subset of $\C^n$ for $n>1$ and $\Lambda\subset \CP^n$ 
is a complex hyperplane such that $\overline E\cap \Lambda=\varnothing$ 
and $\overline E$ is polynomially convex in 
$\CP^n\setminus\Lambda \cong\C^n$, then $\C^n\setminus E$ is Oka.
\end{theorem}

Sections 4 and 5 of the survey \cite{Forstneric2023Indag} 
contains many further examples of Oka domains in Euclidean and 
projective spaces. In particular, there are compact sets in $\C^n$, $n>1$,  
of the form $L=K\cup C$, where $K$ is polynomially convex and 
$C$ is a finite union of rectifiable curves,  
such that $L$ fails to be polynomially convex yet  
$\C^n\setminus L$ is an Oka domain 
(see \cite[Theorem 4.12]{Forstneric2023Indag}).  
The set $L$ in these examples is rationally convex, and  
I do not know examples of non-rationally convex compact sets 
in $\C^n$ with Oka complements.

It seems that all known examples of proper 
Oka domains in Stein manifolds are weakly pseudoconcave.
Since Stein domains in Stein manifolds are exactly the pseudoconvex ones and Oka manifolds are dual to Stein manifolds in the sense explained in the Introduction, this seems natural. However, few explicit results are known.

%
%
\begin{problem}\label{prob:geometric}
Let $K$ be a closed subset of $\C^n$, $n>1$, with smooth boundary $bK$.
\vspace{-1mm} 
\begin{enumerate}[\rm (a)]
\item (The inverse Levi problem.) 
Assuming that $\C^n\setminus K$ is Oka, is $bK$ 
necessarily pseudoconvex?
\vspace{-1mm}
\item Assuming that $K$ is compact and $bK$ is pseudoconvex, is
$\C^n\setminus K$ an Oka domain?
\vspace{-1mm}
\item 
Is every strongly pseudoconcave domain in $\C^n$ 
of the form \eqref{eq:model} an Oka domain?
\end{enumerate}
\end{problem}

An affirmative answer to (a) is known for $n=2$, see
\cite[p.\ 389]{Forstneric2023Indag}.
Model examples of non-pseudoconvex domains are Hartogs figures.
For the closed Hartogs figure 
\[
	H=\big\{(z_1,z')\in \C^n :
	|z_1|\le \delta,\ |z'|	\le 1\big\} 
	\cup \big\{(z_1,z'):|z_1|\le 1,\ 1-\delta\le |z'|\le 1\big\}
\]
with $0<\delta<1$, the complement $\C^n\setminus H$ fails to be Oka
\cite[Example 4.23]{Forstneric2023Indag}. It is not known 
whether the complement of the Hartogs triangle
$
	\{(z_1,z_2)\in\C^2: 0\le |z_1| \le |z_2|\le 1\}
$
is an Oka domain \cite[Problem 4.24]{Forstneric2023Indag}.

In a related direction, Forstneri\v c and Kusakabe \cite{ForstnericKusakabe2025}
studied Oka properties of disc tubes in certain hermitian holomorphic 
line bundles on compact complex manifolds. The following is their main result.

%
%
\begin{theorem}\label{th:density2} 
Let $E$ be a holomorphic line bundle on a compact complex manifold $X$.
Assume that for each point $x\in X$ there exists a divisor $D$ in the
complete linear system $|E|$ 
whose complement $X\setminus D$ is a Stein 
neighbourhood of $x$ with the density property.
Given a semipositive hermitian metric $h$ on $E$, the disc bundle 
$\{e\in E:|e|_h<1\}$ is a pseudoconcave Oka domain while 
$\{e\in E:|e|_h>1\}$ is Kobayashi hyperbolic.
In particular, the zero section of $E$ admits a basis of 
Oka neighbourhoods $\{|e|_h<c\}$ with $c>0$. 
This holds in particular if $E$ is an ample line bundle on 
a rational homogeneous manifold $X$ of dimension $>1$.
\end{theorem}

When $X$ is a complex Grassmannian, 
an ample line bundle $E$ on $X$ satisfies the conditions in the theorem 
with divisors $D\subset X$ such that $X\setminus D\cong \C^n$ 
with $n=\dim X>1$. If $h$ is a semipositive hermitian 
metric on $E$ then the intersection of the unit disc bundle 
$\{|e|_h<1\}$ with any such chart is a pseudoconcave Hartogs
domain of the form
$\Omega=\{(z,t)\in \C^n\times\C: |t|<\phi(z)\}$,  
where $\phi$ is a positive continuous function on $\C^n$
such that $\log \phi$ is plurisubharmonic, and there is a constant $c>0$ such 
that $\phi(z)\ge c\,|z|$ holds for all $z\in\C^n$. By using 
Theorem \ref{th:unboundedsets}
one can show that $\Omega$ is an Oka domain,
and the proof of Theorem \ref{th:density2} is then completed by applying
the localisation Theorem \ref{th:localisation}. Similar arguments
apply when $X$ is a rational homogeneous manifold.

%
%
%
%
\subsection{Oka-1 manifolds.} \label{ss:Oka1}

A complex manifold $Y$ is an {\em Oka-1 manifold} 
if it enjoys the Oka property with approximation and jet interpolation
for maps $X\to Y$ from any open Riemann surface $X$.
In other words, $Y$ contains plenty of holomorphic curves parametrised 
by any open Riemann surface. 
This class of manifolds was introduced and studied 
by Alarc\'on and the author in \cite{AlarconForstneric2025MZ}. 
The following is a special case of \cite[Theorem 2.2]{AlarconForstneric2025MZ}.

\begin{theorem}\label{th:Oka1}
If $Y$ is complex manifold of dimension $n$ and $E\subset Y$
is a closed subset with zero Hausdorff $(2n-1)$-measure
such that $Y$ is dominable by $\C^n$
at every point $y\in Y\setminus E$, then $Y$ is an Oka-1 manifold.
\end{theorem}


A complex manifold $Y$ satisfying the condition in the above theorem is 
said to be {\em densely dominable} by $\C^n$; it is 
{\em strongly dominable} if the condition holds with $E=\varnothing$. 
It is not known whether dense or strong dominability implies
the Oka property for maps from higher dimensional
Stein manifolds. There are examples of noncompact 
Oka-1 manifolds which fail to be Oka, but we do not know any such 
example among compact manifolds.
In \cite{AlarconForstneric2025MZ} we found many examples of Oka-1 manifolds
among compact complex surfaces. 
In particular, every Kummer surface and 
every elliptic K3 surface is Oka-1 (see \cite[Sect.\ 8]{AlarconForstneric2025MZ}). 
It is not known whether any or all of these surfaces are Oka manifolds.
(A survey of known Oka manifolds among compact complex
surfaces can be found in 
\cite[Sect.\ 7.3]{Forstneric2017E}.) 
It is conjectured that every projective rationally
connected manifold is Oka-1
\cite[Conjecture 9.1]{AlarconForstneric2025MZ}.
By Campana and Winkelmann \cite{CampanaWinkelmann2023},
every such manifold admits holomorphic lines 
$\C\to X$ with dense images.

New examples and functorial properties of Oka-1 manifolds
were presented by L\'arusson and the author in 
\cite{ForstnericLarusson2025MZ}.
We also formulated and studied the algebraic version of the Oka-1 condition, 
called aOka-1, which concerns approximation and 
interpolation of holomorphic maps from open subsets 
of affine algebraic curves to algebraic manifolds by regular algebraic maps.
Every aOka-1 manifold is also an Oka-1 manifold.
The aOka-1 property is a birational invariant for compact 
algebraic manifolds, and it holds for all rational manifolds
and all algebraically elliptic projective manifolds 
(see \cite[Theorem 1.6]{ForstnericLarusson2025MZ}). 

It is easily seen that every projective aOka-1 manifold is rationally connected 
\cite[Proposition 1.7]{ForstnericLarusson2025MZ}.
Conversely, Benoist and Wittenberg
\cite[Theorem 1.2]{BenoistWittenberg2025} proved that every 
projective {\em rationally simply connected manifold} is an aOka-1 manifold, 
and this class contains every smooth
hypersurface of degree $d$ in $\CP^n$ with ${n\geq d^2-1}$
(see \cite[Corollary 1.3]{BenoistWittenberg2025}).
Note that their {\em tight approximation property} 
is equivalent to the aOka-1 property.

%
%
%
%
\subsection{The Oka principle on families of open Riemann surfaces.}\label{ss:families}
In \cite{Forstneric2024Runge}, the Oka principle was established for maps from very general families of open Riemann surfaces to Oka manifolds. 
Assume that $X$ is a smooth open orientable surface,
$B$ is a finite CW complex or a countable 
locally compact CW-complex of finite dimension, and 
$\{J_b\}_{b\in B}$ is a continuous family of complex structures on 
$X$ of local H\"older class $\Cscr^\alpha$, $0<\alpha<1$.
A continuous map $f:B\times X\to Y$ to a complex manifold 
$Y$ is said to be fibrewise holomorphic if the
map $f(b,\cdotp):X\to Y$ is holomorphic from the open 
Riemann surface $(X,J_b)$ for every $b\in B$. The following 
Oka principle is a special case of \cite[Theorem 1.6]{Forstneric2024Runge}.

%
%
\begin{theorem}\label{th:families}
(Assumptions as above.) If $Y$ is an Oka manifold 
then every continuous map $f:B\times X\to Y$
is homotopic to a continuous fibrewise holomorphic map $F:B\times X\to Y$. 
If in addition $B$ is a manifold of class $\Cscr^l$
and the family of complex structures $\{J_b\}_{b\in B}$
is of local class $\Cscr^{l,(k,\alpha)}(X)$,
where $0\le l\le k+1$ and $0<\alpha<1$, then $F$ can be chosen of local class
$\Cscr^{l,(k+1,\alpha)}$.  
\end{theorem}

The statement concerning the map $F$ 
in the second part of the theorem means that 
$F$ has $l$ derivatives in the parameter $b\in B$ followed by $k+1$ 
derivatives in the space variable $x\in X$, and these
derivatives are of local H\"older class $\Cscr^\alpha$ in $x\in X$.
The analogous definition applies to the family $\{J_b\}_{b\in B}$.

The more precise result \cite[Theorem 1.6]{Forstneric2024Runge}
includes approximation of fibrewise holomorphic maps
on suitable families of compact Runge subsets of $X$.
Theorem \ref{th:families} relates the Oka theory to Teichm\"uller spaces.
It was also used in the construction of 
families of immersed conformal minimal surfaces 
$X\to \R^n$ and holomorphic null curves $X\to\C^n$, $n\ge 3$,
with a prescribed family of conformal structures on the surface $X$;
see \cite[Sect.\ 8]{Forstneric2024Runge}. 
Further applications are expected.

The proof of Theorem \ref{th:families} proceeds by inductively extending
the domain in $B\times X$ on which a given map is fibrewise holomorphic,
using approximation on closed subsets with compact
Runge fibres to ensure convergence to a limit map.
It also uses solutions of the Beltrami equation on an open Riemann surface
$(X,J)$ to show that a small variation of the complex structure 
$J$ on a smoothly bounded domain $\Omega\Subset X$ 
can be realised by a small variation 
of $\Omega$ in $X$ with its original structure $J$. This is 
special case of Hamilton's theorem \cite{Hamilton1977} but
with more precise regularity of the maps. 
Finally, the problem is reduced locally in the parameter $b\in B$ 
and semiglobally in the space variable $x\in X$ to the Oka principle
given by Theorem \ref{th:Okamain}. 

An important fact used 
in the proof of Theorem \ref{th:families} is that the condition 
on a compact set $K$ in a surface $X$ to be Runge is purely topological
(its complement $X\setminus K$ has no relatively compact connected
component), hence independent of the choice of the complex 
structure on $X$. This is no longer the case if $X$
is a smooth open manifold of dimension $2n>2$ endowed 
with a family of Stein structures $\Jscr=\{J_b\}_{b\in B}$. 
The following example shows that one must impose a suitable 
condition on the family $\Jscr$ to get an analogue of 
Theorem \ref{th:families}. The proof will appear in a subsequent
publication (in preparation).

%
%
\begin{theorem}\label{th:nontame}
Given a compact set $K\subset \R^{2n}$ $(n>1)$ with nonempty interior,
there is a family of smooth integrable complex structures 
$\{J_t\}_{t\in\R}$ on $\R^{2n}$, depending smoothly on $t\in \R$, 
such that $J_0$ is the standard complex structure on $\C^n$,
the manifold $(\R^{2n},J_t)$ is biholomorphic to $\C^n$ for every $t\in\R$, 
and for any neighbourhood $U\subset \R$ of $0\in\R$ the 
set $\bigcup_{t\in U} \wh{K}_{J_t}\subset \R^{2n}$ is unbounded. 
\end{theorem}

This phenomenon excludes the possibility 
of any reasonable analysis on such a family of Stein structures. 
The condition that one must impose to obtain positive results
is that for any compact set $K$ in $X$, the $J_b$-convex hull
$\wh K_{J_b}\subset X$ \eqref{eq:hull} is an upper semicontinuous 
set-valued function of the parameter $b\in B$. Equivalently,
there is a continuous function $\rho:B\times X\to [0,\infty)$ such 
that for every $b\in B$, the function $\rho(b,\cdotp):X\to [0,\infty)$ 
is a strongly $J_b$-plurisubharmonic exhastion. 
A paper on this topic is in preparation.

%
%
%
%
\section{Elliptic and subelliptic manifolds.}\label{sec:elliptic}

Gromov introduced the following conditions in \cite[0.5, p.\ 855]{Gromov1989}.

%
%
\begin{definition}\label{def:elliptic}
A complex manifold $Y$ is elliptic if there is a holomorphic
vector bundle $\pi:E\to Y$ and a holomorphic map $s:E\to Y$
such that $s$ restricts to the identity map on the zero section 
$E(0)$ of $E$, and for any point $y\in Y$ the differential
$ds_{0_y}:T_{0_y}E\to T_yY$ maps the fibre 
$E_y=\pi^{-1}(y)$ onto $T_yY$. An algebraic manifold $Y$
is algebraically elliptic if the vector bundle $\pi:E\to Y$ and
the map $s:E\to Y$ can be chosen algebraic.
\end{definition}

Here, $0_y$ denotes the origin of the vector space $E_y=\pi^{-1}(y)$.
A map $s$ with these properties is called a {\em dominating spray} on $Y$.
An ostensibly weaker condition, {\em subellipticity} 
(see \cite[Def.\ 2]{Forstneric2002MZ}),
asks for the existence of a finite family of holomorphic
sprays $(E_j,\pi_j,s_j)$ on $Y$ $(j=1,\ldots,m)$ 
which is dominating in the sense that
\begin{equation}
\label{eqn:domination}
    (ds_1)_{0_y}(E_{1,y}) + (ds_2)_{0_y}(E_{2,y})+\cdots
                     + (ds_m)_{0_y}(E_{m,y})= T_y Y
                     \quad \text{for all $y\in Y$}.
\end{equation}
An algebraic manifold $Y$ is said to be {\em algebraically subelliptic} 
if the same condition holds with algebraic sprays. 
See also \cite[Def.\ 5.6.13]{Forstneric2017E}.
Examples of elliptic and subelliptic manifolds can be found 
in \cite{Gromov1989}, \cite[Sect.\ 6.4]{Forstneric2017E}, and in the surveys \cite{Forstneric2023Indag,ForstnericLarusson2011}.
We recall a few.

%
%
\begin{example} \label{example:elliptic}
\noindent (A) 
{\em Every complex homogeneous manifold is elliptic.}
Assume that a complex Lie group $G$ acts on a complex 
manifold $Y$ transitively by holomorphic automorphisms.
Let $\ggot\cong\C^p$ denote the Lie algebra of $G$ and
$\exp:\ggot\to G$ the exponential map.
The holomorphic map $s:Y\! \times\ggot \cong Y\times\C^p \to Y$ 
given by
$ 
      s(y,v)= \exp v\, \cdotp y \in Y (y\in Y,\ v\in \ggot)
$ 
is a dominating spray on $Y$. Such sprays were used by 
Grauert \cite{Grauert1957I,Grauert1957II,Grauert1958MA}.

\smallskip
\noindent (B) 
{\em If $Y$ is a complex manifold whose tangent bundle is 
spanned by finitely many $\C$-complete holomorphic vector fields
$V_1,\ldots,  V_m$, then $Y$ is elliptic.} 
Indeed, let $\phi_j^t$, $t\in \C$, be the flow of $V_j$. 
The map $s:Y\times \C^m \to Y$ given by
\begin{equation}\label{eq:flow-spray}
      s(y,t)= s(y,t_1,\ldots,t_m) =
      \phi_1^{t_1}\circ \phi_2^{t_2}\circ\cdots\circ \phi_m^{t_m}(y)
\end{equation}
clearly satisfies $s(y,0)=y$ and $\frac{\di}{\di t_j} s(y,0)= V_j(y)$
for any $y\in Y$ and $j=1,\ldots,m$. Thus, $s$ is dominating 
precisely when the vectors $V_1(y),\ldots,  V_m(y)$ 
span the tangent space $T_y Y$ for every $y\in Y$. 

\smallskip
\noindent (C) Assume that $\pi:E\to Y$ is a holomorphic vector bundle
and $\phi:E\to E$ is a holomorphic map which restricts to the
identity on the zero section $E(0)$ of $E$. Then, the map
$s=\pi\circ\phi:E\to Y$ is a spray on $Y$. Let us identify $Y$ 
with $E(0)$ and the restricted tangent bundle
$TE|_{Y}$ with $TY\oplus E$. The spray $s$ is dominating
at $y\in Y$ if for every $v\in T_yY$ there is a
point $e \in E_y$ such that $ds_y(e) = (d\pi)_y \circ (d\phi)_y (e)=v$. 
\end{example}

The importance of ellipticity and subellipticity lies in the following result.

%
%
\begin{theorem}\label{th:ellCAP}
\begin{enumerate} [\rm (a)] 
\item {\rm (Gromov \cite[0.6, p.\ 855]{Gromov1989})} 
Every elliptic manifold is an Oka manifold. 
\vspace{-1mm}
\item {\rm (\cite[Theorem 1.1]{Forstneric2002MZ})} 
Every subelliptic manifold is an Oka manifold. 
\end{enumerate}
\end{theorem}

The modern proof of Theorem \ref{th:ellCAP} 
(see \cite[Chap.\ 6]{Forstneric2017E}) consists of two parts. 
The first part is the implication
\[
	\text{(sub)elliptic} \Longrightarrow \text{h-Runge approximation}
	\Longrightarrow \text{CAP}. 
\]
Here, a complex manifold $Y$ is said to satisfy the h-Runge approximation condition if the following holds. Given a pair $K\subset L$ of compact 
$\Oscr(X)$-convex sets in a Stein space $X$ and 
a homotopy of holomorphic maps $f_t:U\to Y$, $t\in I=[0,1]$, 
of holomorphic maps from a neighbourhood of $K$
such that $f_0$ extends to a holomorphic
map from a neighbourhood of $L$ in $Y$, 
we can approximate the homotopy $f_t$ uniformly on $K$ by a 
homotopy of holomorphic maps $\tilde f_t:L \to Y$, $t\in I$, 
such that $\tilde f_0=f_0$. This is a fairly elementary consequence of
subellipticity and the parametric Oka--Weil theorem for sections
of vector bundles; see \cite[Theorem 6.6.1]{Forstneric2017E}. 
Since a compact convex set $K\subset\C^n$ is holomorphically
contractible, every holomorphic map $f$ from a convex neighbourhood
of $K$ is homotopic to a constant map through holomorphic maps.
Hence, the h-Runge approximation property of $Y$ implies
that $f$ is a limit of entire maps $\C^n\to Y$, 
which means that $Y$ enjoys CAP.

The implication $\text{CAP} \Longrightarrow \text{OKA}$
is highly nontrivial; see \cite[Chap.\ 5]{Forstneric2017E} for the proof.

The notion of (sub)ellipticity is also defined for holomorphic
submersions $h:Z\to X$, and it implies the parametric Oka property
for sections $X\to Z$ provided that $X$ is Stein
(see \cite[Sect.\ 4]{Gromov1989}, \cite{ForstnericPrezelj2002}, 
and \cite[Chap.\ 6]{Forstneric2017E}). An even more general case
pertains to sections of stratified subelliptic submersions over
Stein spaces; see \cite{Forstneric2010PAMQ} and 
\cite[Theorem 6.2.2]{Forstneric2017E}. There is also
a relative Oka principle for sections of ramified holomorphic 
maps onto Stein spaces; see \cite{Forstneric2003FM} and 
\cite[Sect.\ 6.14]{Forstneric2017E}. 

It is easily seen that every Stein Oka manifold is 
elliptic \cite[3.2.A]{Gromov1989}.
Gromov asked whether every Oka manifold is elliptic. 
The first counterexamples for noncompact manifolds 
were given by Kusakabe in \cite{Kusakabe2020PAMS}. He showed in 
particular that if $S$ is a closed tame countable set in $\C^n$, $n>1$, 
whose set of limit points is discrete then $\C^n\setminus S$ is an
Oka domain (see \cite[Theorem 1.2]{Kusakabe2020PAMS}) but
there are nonelliptic domains of this type. For example,
taking $Z=\{0\}\cup \{1/j:j\in\N\}\subset \C$, the domain
$\C^n\setminus (Z^2 \times \{0\}^{n-2})$ for $n\ge 3$ 
is Oka but not elliptic or subelliptic \cite[Corollary 1.4]{Kusakabe2020PAMS}.  
Theorem \ref{th:complements}, together 
with results of Andrist, Shcherbina and Wold \cite{AndristShcherbinaWold2016},  
gives a much more abundant class of examples in the following
corollary. 

%
%
\begin{corollary}\label{cor:complements}
Let $n\ge 3$. For every compact polynomially convex set $K\subset \C^n$
with infinitely many limit points the complement $\C^n\setminus K$ is 
Oka but not subelliptic. The analogous result holds for compact 
holomorphically convex sets in any Stein manifold 
with the density property of dimension $\ge 3$.
\end{corollary}

The following recent result of Forstneri\v c and L\'arusson
\cite{ForstnericLarusson2025MRL} confirms 
Gromov's conjecture \cite[3.2.A']{Gromov1989}. 

%
%
\begin{theorem}\label{th:projective}
Every projective Oka manifold is elliptic.
\end{theorem}

From Theorems \ref{th:ellCAP} (b) and \ref{th:projective} we conclude
the following.

%
%
\begin{corollary}\label{cor:projective}
For a projective manifold, being Oka, elliptic, or subelliptic are 
equivalent conditions.
\end{corollary}

There are projective Oka manifolds which fail to be
algebraically elliptic; for example, abelian varieties. 
Hence, the algebraic counterpart to Theorem \ref{th:projective} 
is not true and Serre's GAGA principle fails in this problem.

The spray bundle $E\to Y$ in the proof of 
Theorem \ref{th:projective} is the direct sum of copies 
of the dual of a sufficiently ample line bundle on the projective 
manifold $Y$ (hence, $E$ is Griffiths negative), and the spray
map $s:E\to Y$ is of the type described in
Example \ref{example:elliptic} (C). The construction in 
\cite{ForstnericLarusson2025MRL} first produces a local dominating
spray $s_0:U\to Y$ defined on an open neighbourhood $U\subset E$
of the zero section $E(0)$; this does not use any 
additional hypothesis on $Y$. Since $E$ is negative,
it is a $1$-convex manifold with the exceptional subvariety $E(0)\cong Y$. 
Assuming that $Y$ is an Oka manifold, the existence of a global 
holomorphic map $s:E\to Y$ which agrees with $s_0$ 
to the second order along
$E(0)$ then follows from the Oka principle for maps
from 1-convex manifolds to Oka manifolds, due to 
Stopar \cite{Stopar2013}.
Note that such $s$ is a global dominating spray on $Y$,
thus proving that $Y$ is elliptic.

%
%
Algebraic analogues of Stein manifolds are affine algebraic manifolds. 
The properties of holomorphic functions on Stein spaces 
(see (i)--(iv) in Sect.\ \ref{sec:Stein}) also 
hold for (regular) algebraic functions on affine algebraic varieties,
and they can be used to define algebraic Oka properties of algebraic manifolds. We discuss them briefly, referring to \cite[Sect.\ 6]{Forstneric2023Indag} for a more comprehensive presentation.

Recall that a complex algebraic manifold $Y$ is said to be algebraically
elliptic if it admits an algebraic dominating spray $s:E\to Y$
on an algebraic vector bundle $\pi:E\to Y$. 
Similarly, $Y$ is said to be algebraically subelliptic if it admits 
finitely many algebraic sprays $(E_j,\pi_j,s_j)$ satisfying 
\eqref{eqn:domination}. Gromov proved that algebraic subellipticity 
is a Zariski local property \cite[3.5 B]{Gromov1989}, and 
it is stable with respect to removing algebraic 
subvarieties of codimension $>1$ \cite[3.5 C]{Gromov1989}. 
(Removal of hypersurfaces is not allowed in general as it may lead
to hyperbolic manifolds.)
Examples and properties of such manifolds, with further 
references, can be found in \cite[Sect.\ 6.4]{Forstneric2017E}
and \cite[Sect.\ 6]{Forstneric2023Indag}. 
The following recent result is due to Kaliman and Zaidenberg
\cite{KalimanZaidenberg2024FM}.

%
%
\begin{theorem}\label{th:KZelliptic}
Every algebraically subelliptic manifold is algebraically elliptic. 
\end{theorem}

Hence these two properties are equivalent, so 
algebraic ellipticity is a Zariski-local condition.
No such result is known in the holomorphic category. 
See also \cite[Theorem 6.2]{Forstneric2023Indag}
for some other equivalent properties. 

The main Oka-type property of algebraically elliptic manifolds
is the algebraic homotopy approximation theorem for maps from 
affine algebraic varieties; see \cite[Theorem 3.1]{Forstneric2006AJM}, 
\cite[Theorem 6.15.1]{Forstneric2017E}, and the recent 
generalisations in \cite[Sect.\ 2]{AlarconForstnericLarusson2024}.
In particular, we have the following result.

\begin{theorem} \label{th:a-HAT}
Let $X$ be an affine algebraic variety and $Y$ be an 
algebraically elliptic manifold. Then, every holomorphic 
map $X\to Y$ homotopic to an algebraic map 
is a limit of algebraic maps in the compact-open topology.
\end{theorem}

It follows that every algebraically elliptic manifold $Y$ satisfies 
the algebraic convex approximation property, aCAP: every 
holomorphic map from a convex set $U\subset \C^n$ to $Y$ 
is a limit of algebraic maps $\C^n\to Y$ uniformly on compacts.
Theorem \ref{th:a-HAT} can be used to prove the following 
result on algebraic domination. The first 
part is a special case of \cite[Theorem 1.6]{Forstneric2017Indam},
while the second part is due to Kusakabe 
\cite[Theorem 1.2]{Kusakabe2022surjective}.

\begin{theorem}\label{th:algdomination}
Let $Y$ be an algebraically elliptic manifold of dimension $n$.
If $Y$ is compact then it admits an algebraic map $f:\C^n\to Y$
satisfying $f(\C^{n} \setminus \mathrm{br}(f))=Y$, where 
$\mathrm{br}(f)$ denotes the branch locus of $f$.
If $Y$ is noncompact then it admits an algebraic map $\C^{n+1}\to Y$ 
satisfying $f(\C^{n+1} \setminus \mathrm{br}(f))=Y$.
\end{theorem}

There are examples of holomorphic maps  
from an affine algebraic manifold $X$ to an  
algebraically elliptic manifold $Y$ which are not homotopic 
to any algebraic map
(see \cite[Examples 6.15.7, 6.15.8]{Forstneric2017E}). 
Thus, the algebraic basic Oka property, aBOP, 
fails in these examples. The following result 
of L\'arusson and Truong \cite[Theorem 2]{LarussonTruong2019}
shows that this a rather common phenomenon, 
and Theorem \ref{th:a-HAT} might be the closest analogue 
of the Oka principle in the algebraic category.  

%
%
\begin{theorem}\label{th:P1isbad}
If $Y$ is an algebraic manifold which contains a rational curve  
or is compact, then $Y$ does not have the algebraic basic Oka property,
aBOP.
\end{theorem}

I am not aware of a single example of an algebraic 
manifold $Y$ with nontrivial topology satisfying aBOP. 
The simplest example to consider is $\C^2\setminus\{0\}$,
which is algebraically elliptic.

By \cite[Theorem 3]{LarussonTruong2019} of L\'arusson and Truong,
together with Theorem \ref{th:KZelliptic},
every smooth nondegenerate toric variety is algebraically elliptic.
A result of Banecki \cite{Banecki2024X} says that 
every rational projective manifold is algebraically elliptic.
This generalises the result of Arzhantsev et al. 
\cite[Theorem 1.3]{ArzhantsevKalimanZaidenberg2024}
that every uniformly rational projective manifold is algebraically
elliptic. See also the recent examples 
\cite{KalimanZaidenberg2024MRL,KalimanZaidenberg2024X}
and the survey \cite{Zaidenberg2024ellipticitysurvey}.
A major source of examples of algebraically elliptic manifolds are 
{\em flexible manifolds}; see 
\cite{ArzhantsevFlennerKalimanKutzschebauchZaidenberg2013DMJ}.
These are algebraic manifolds whose tangent space at every point is 
spanned by algebraic vector fields with complete algebraic flows,
so they admit dominating algebraic sprays of type \eqref{eq:flow-spray}.
A list of examples of such manifolds can be found in \cite[p.\ 394]{Forstneric2023Indag}.

\section{The h-principles for minimal surfaces in $\R^n$ 
and holomorphic null curves in $\C^n$.}
\label{sec:minimal}

We begin by recalling the basic facts; 
see \cite{AlarconForstnericLopez2021,Osserman1986}.
Let $M$ be a connected open Riemann surface.
A smooth immersion $u=(u_1,\ldots,u_n): M\to \R^n$, $n>1$, is 
conformal (angle preserving) if and only if its $(1,0)$-differential 
$\di u=(\di u_1,\ldots, \di u_n)$ 
(the $\C$-linear part of the differential $du$)
satisfies the nullity condition
\begin{equation}\label{eq:nullity}
	(\di u_1)^2 + (\di u_2)^2 + \cdots + (\di u_n)^2 = 0.
\end{equation}
Let $n\ge 3$.
A conformal immersion $u: M\to \R^n$ is {\em minimal}, in the sense 
that it parametrizes a minimal surface in $\R^n$, if and only if it is 
harmonic, $\dibar\di u=0$;
this holds if and only if $\di u$ is a holomorphic $(1,0)$-form. 
Such an immersion $u$ is said to be {\em nonflat} if
the image $u(M)\subset\R^n$ is not contained in an affine $2$-plane. 
We denote by $\Mgot(M,\R^n)$ the space of conformal minimal 
immersions $M\to\R^n$ with the compact-open topology, and by 
$\Mgot_{\mathrm{nf}}(M,\R^n)$ the subspace of $\Mgot(M,\R^n)$ 
consisting of all nonflat conformal minimal immersions. 

A holomorphic immersion $F: M\to \C^n$, $n\ge 3$, is 
said to be a {\em null curve}
if the differential $dF=\di F=(dF_1,\ldots, dF_1)$ satisfies 
the nullity condition \eqref{eq:nullity}.
Such $F$ is said to be {\em nonflat} if $F(M)$ is not contained in 
an affine complex line. Denote by $\Ngot(M,\C^n)$ the spaces
of null holomorphic immersions $M\to\C^n$ with the compact-open
topology, and by $\Ngot_{\mathrm{nf}}(M,\C^n)$ 
the subspace consisting 
of all nonflat null holomorphic immersions $M\to\C^n$. 
Since $dF=2\di (\Re F)$, the real part $u=\Re F$ 
of a (nonflat) null curve is a (nonflat) conformal minimal immersion 
$M\to \R^n$; the converse holds if $M$ is simply connected. 
Hence, we have natural inclusions 
\[
	\Re \Ngot(M,\C^n) \hra \Mgot(M,\R^n),\qquad 
	\Re \Ngot_{\mathrm{nf}}(M,\C^n) \hra \Mgot_{\mathrm{nf}}(M,\R^n), 
\]
where $\Re \Fscr$ denotes the space of real parts of maps in 
a space $\Fscr$. 

Choose a nowhere vanishing holomorphic 
$1$-form $\theta$ on $M$ (such exists by the Oka--Grauert principle 
\cite[Theorem 5.3.1]{Forstneric2017E}). 
Given a conformal minimal immersion $u: M\to \R^n$, the map
$
	f=2\di u/\theta =(f_1,\ldots, f_n): M\to \C^n
$ 
is holomorphic since $u$ is harmonic, and by \eqref{eq:nullity} 
it has range in $\A_*=\A\setminus\{0\}$ where $\A$ is the null quadric 
\begin{equation}\label{eq:nullq}
	\A = \{(z_1,\ldots,z_n) \in\C^n : 
	z_1^2+z_2^2+\cdots + z_n^2=0\}.
\end{equation}
Conversely, a holomorphic map $f: M\to \A_*$ such that the 
1-form $f\theta$ has vanishing real periods (i.e., 
$\int_\gamma \Re(f\theta)=0$ holds for every closed curve $\gamma$ in $M$)
determines a conformal minimal immersion $u: M\to\R^n$ by 
$u(x)= \int^x \Re(f\theta)$, $x\in M$. If $f\theta$ has vanishing 
complex periods then it integrates to a holomorphic null curve 
$F(x)=\int^x f\theta$. This is the classical 
{\em Weierstrass representation} of minimal surfaces and null curves;
see \cite[Sec.\ 2.3]{AlarconForstnericLopez2021} or \cite{Osserman1986}.

Consider the following commuting diagram where 
the first vertical arrow is the real part projection, the map 
$\vartheta$ is given by $F\mapsto \di F/\theta$, and 
$\psi$ is given by $u\mapsto 2\di u/\theta$.
\begin{equation}\label{eq:diagram}
\xymatrix{
	\Ngot_{\mathrm{nf}}(M,\C^n)  \ar[r]^\vartheta \ar[d]  
	&  \Oscr(M,\A_*) \ \ar@{^{(}->}[r]  & \Cscr(M,\A_*)  \\ 
	\Re\Ngot_{\mathrm{nf}}(M,\C^n) \  \ar@{^{(}->}[r]   
	&  \Mgot_{\mathrm{nf}}(M,\R^n)  \ar[u]_\psi & 
}
\end{equation}
The punctured null quadric $\A_*=\A\setminus\{0\}$ 
is a homogeneous space of the complex Lie group $\C^*\times O(n,\C)$, 
where $O(n,\C)=\{A\in GL(n,\C): AA^t=I\}$ 
is the orthogonal group over $\C$.
Hence, the inclusion $\Oscr(M,\A_*) \longhookrightarrow \Cscr(M,\A_*)$ 
of the space of holomorphic maps $M\to \A_*$ in 
the space of continuous maps satisfies the Oka principle
and so is a weak homotopy equivalence. 
By using elements of the convex integration theory, 
one shows that every holomorphic map $f:M\to \A_*$ can be deformed
to a map for which the 1-form $f\theta$ is exact holomorphic, 
so it integrates to a null holomorphic immersion 
$\int f\theta:M\to \C^n$. Ensuring only that $\Re(f\theta)$ 
is exact yields a conformally immersed minimal surface
$\int \Re (f\theta) :M\to \R^n$. Furthermore, such deformations can be
made in families. This leads to the following result of 
Forstneri\v c and L\'arusson 
\cite[Theorems 1.1 and 4.1]{ForstnericLarusson2019CAG}.

%
%
\begin{theorem}\label{th:WHE1}
All maps in the diagram \eqref{eq:diagram} satisfy the parametric h-principle with 
approximation, and hence they are weak homotopy equivalences.
\end{theorem}

The basic h-principles with approximation were proved in this context by 
Alarc\'on and Forstneri\v c in \cite{AlarconForstneric2014IM} 
for the maps $\vartheta$ and $\psi$, and in  
\cite[Theorem 1.1]{AlarconForstneric2018Crelle} for the inclusion 
$\Re\Ngot_{\mathrm{nf}}(M,\C^n) \hra \Mgot_{\mathrm{nf}}(M,\R^n)$.
It is not know whether the inclusion 
$\Re \Ngot(M,\C^n) \hra \Mgot(M,\R^n)$
satisfies the h-principle. 

Let $ds^2$ denote the Euclidean metric on $\R^n$.
A minimal surface $u:M\to\R^n$ is said to be {\em complete} if the Riemannian 
metric $g=u^* ds^2$ induces a complete metric on $M$.
Completeness is an important and much studied property in the global theory 
of minimal surfaces. Alarc\'on and L\'arusson \cite{AlarconLarusson2025JGEA}
proved that the inclusion 
$\Mgot^c_{\mathrm{nf}}(M,\R^n) \hra \Mgot_{\mathrm{nf}}(M,\R^n)$ 
of the space of {\em complete} nonflat minimal immersions $M\to\R^n$ 
in $\Mgot_{\mathrm{nf}}(M,\R^n)$ also satisfies the parametric homotopy principle,
whence is a weak homotopy equivalence. When $M$ is of finite topological type, 
the inclusion is a genuine homotopy equivalence. 
The following corollary 
(see \cite[Corollary 1.6]{ForstnericLarusson2019CAG}) is 
seen by analysing the path components of the space $\Cscr(M,\A_*)$,
a purely topological issue. 

%
%
%
%
\begin{corollary} \label{cor:components}
Let $M$ be a connected open Riemann surface with 
$H_1(M;\Z)\cong \Z^l$, $l\in \Z_+\cup\{\infty\}$. 
Then the path connected components of each of the spaces 
$\Mgot_{\mathrm{nf}}(M,\R^3)$, $\Mgot^c_{\mathrm{nf}}(M,\R^3)$, 
and $\Ngot_{\mathrm{nf}}(M,\C^3)$ 
are in one-to-one correspondence with the elements 
of the abelian group $(\Z_2)^l$. 
If $n\ge 4$ then these spaces are path connected.
\end{corollary} 

The {\em total curvature} of a conformal minimal surface $u:M\to\R^n$ is the 
integral $\TC(u) = \int_M K\, d\sigma$ of the Gaussian curvature function 
$K:M\to (-\infty,0]$ of $u$ with respect to the area measure $d\sigma$ 
on $M$ induced by the immersion $u$. 
A minimal surface $u:M\to\R^n$ is said to have
{\em finite total curvature} if $\TC(u)>-\infty$. 
Complete minimal surfaces of finite total Gaussian curvature 
are among the most intensively studied minimal surfaces
which play an important role in the classical global theory;
see \cite{Osserman1986} and \cite[Chap.\ 4]{AlarconForstnericLopez2021},
among other sources.
In fact, most classical examples of minimal surfaces are of this kind. 
Although this family of minimal surfaces has been a focus of interest 
since the seminal work of Osserman in the 1960s, 
the theories of approximation and interpolation 
for complete minimal surfaces of finite total curvature in $\mathbb{R}^n$
have been developed only recently; see \cite{AlarconCastroInfantesLopez2019CVPDE,AlarconLopez2022APDE,Lopez2014TAMS}.  
The paper \cite{AlarconForstnericLarusson2024} 
by Alarc\'on, Forstneri\v c, and L\'arusson 
provides the first contributions to the homotopy theory of this 
important class of minimal surfaces, which we now present.

The null quadric $\A$ \eqref{eq:nullq} defines a holomorphic  
subbundle $\mathcal A$ of the vector bundle 
$(T^*M)^{\oplus n}$ with fibre $\A_*$ 
whose sections are $n$-tuples $\phi=(\phi_1,\ldots,\phi_n)$ 
of $(1,0)$-forms on $M$ without common zeros such that the 
map $[\phi_1:\cdots:\phi_n]:M\to \CP^{n-1}$ 
takes values in the projective quadric 
\begin{equation} \label{eq:projquadric}
	Q = \bigl\{ [z_1:z_2:\cdots : z_n]\in \CP^{n-1}: 
		z_1^2+z_2^2+ \cdots + z_n^2=0\bigr\}.
\end{equation}
If $u:M\to\R^n$ is a conformally immersed minimal surface 
with finite total curvature, then $M$ is the complement in a compact 
Riemann surface $\overline M$ of a nonempty finite set 
$E=\{x_1, \ldots, x_m\}$ whose points are called the ends of $M$
(such $M$ is called an {\em affine Riemann surface}), the bundle
$\Acal\to M$ is algebraic (it is holomorphically trivial 
but not necessarily algebraically trivial), 
and $\phi=\di u$ is a meromorphic $1$-form on $\overline M$ 
without zeros or poles on $M$, that is, an algebraic section 
of $\Acal$ over $M$. Furthermore, the immersion $u$ is complete 
if and only if $\di u$ has an effective pole at every point of 
$E=\overline M\setminus M$. 
The surface $u(M)$ is then properly immersed in $\R^n$ 
and has a fairly simple asymptotic behaviour at every end
of $M$, described by the Jorge--Meeks theorem \cite{JorgeMeeks1983T}.
Denote by $\Ascr^1(M,\A_*)$ the space of meromorphic 
$1$-forms $\phi=(\phi_1,\ldots,\phi_n)$ on $\overline M$ having no zeros or 
poles in $M$ and satisfying \eqref{eq:nullity}, 
that is, algebraic sections of the bundle $\mathcal A\to M$.
Consider the diagram
\begin{equation}\label{eq:diagram2}
\xymatrix{
	\Re \Ngot_*(M,\C^n)  \ \ar@{^{(}->}[r] \ar[dr]_\di  
	&  \Mgot_*(M,\R^n) \ar[d]^\di  \\
	&   \Ascr^1(M,\A_*)
}
\end{equation}
where $\di$ is the $(1,0)$-differential, $\Ngot_*(M,\C^n)$
is the space of nonflat complete (proper) algebraic null immersions
$M\to \C^n$, and $\Mgot_*(M,\R^n)$ is the space of nonflat 
conformal minimal immersions of finite total curvature.
These spaces are endowed with the compact-open topology. 
The following is \cite[Theorem 1.1]{AlarconForstnericLarusson2024}.

%
%
\begin{theorem}\label{th:whe2}
If $M$ is an affine Riemann surface then the maps in \eqref{eq:diagram2}
are weak homotopy equivalences. 
\end{theorem}

The proof in \cite{AlarconForstnericLarusson2024} 
builds upon that of Theorem \ref{th:WHE1}, using 
a multiparameter extension of the algebraic approximation
theorem in Theorem \ref{th:a-HAT}, together with 
a method invented by Alarc\'on and L\'arusson 
\cite{AlarconLarusson2025Crelle} to ensure 
the existence of effective poles at the ends of $M$ (to get completeness). 
For the inclusion $\Re \Ngot_*(M,\C^n) \hookrightarrow \Mgot_*(M,\R^n)$,
Theorem \ref{th:whe2} says in particular that every complete nonflat 
conformal minimal immersion $M\to\R^n$ of finite total curvature 
can be deformed through maps of the same type to the real part of 
a proper algebraic null curve $M\to\C^n$. In addition, one can control 
the flux along an isotopy in $\Mgot_*(M,\R^n)$; see 
\cite[Theorem 1.2]{AlarconForstnericLarusson2024}. 


The Gauss map of a conformal minimal immersion  $u=(u_1,u_2,\ldots,u_n):M\to\R^n$ is the holomorphic map
\[
	G(u)= [\di u_1 : \di u_2 : \cdots : \di u_n]: M\to Q\subset \CP^{n-1}
\] 
in the projective hyperquadric \eqref{eq:projquadric}.
When $n=3$, $Q$ is a quadratic rational curve in $\CP^2$, hence
biholomorphic to the Riemann sphere $\CP^1$. In this case,
the Gauss map is identified with a meromorphic function on $M$,
and it corresponds to the standard Gauss map 
when identifying $\CP^1$ with the 2-sphere.
See \cite[Sect.\ 2.5]{AlarconForstnericLopez2021} for the details. 
The following result of Alarc\'on, Forstneri\v c, and L\'opez \cite{AlarconForstnericLopez2019JGEA}
(see also \cite[Theorem 5.4.1]{AlarconForstnericLopez2021})
can be seen as an h-principle for the Gauss map of minimal surfaces.

%
%
\begin{theorem}\label{th:Gaussmap}
Let $M$ be an open Riemann surface and $n\geq 3$ be an integer.
For every holomorphic map $\Gscr : M\to Q\subset \CP^{n-1}$ into 
the quadric \eqref{eq:projquadric} there exists a conformal minimal immersion 
$u : M\to\R^n$ with the Gauss map $G(u)=\Gscr$ 
and vanishing flux (i.e., $u$ is the real part of a holomorphic null curve
$M\to\C^n$). If $\Gscr(M)$ is not contained in any proper projective subspace
of $\CP^{n-1}$ then $u$ can be chosen to have arbitrary flux and 
to be an embedding if $n\ge 5$ and an immersion with simple 
double points if $n=4$. 
\end{theorem}

The theory of the Gauss map of minimal 
surfaces has a long and rich history, and it came as a surprise that 
every natural candidate map is the Gauss map of some minimal surface.
Further results on the space of Gauss maps of complete minimal surfaces
were obtained by Alarc\'on and L\'arusson \cite{AlarconLarusson2024complete}.



\end{document}